# ESTIMATING DEFORMATIONS OF ISOTROPIC GAUSSIAN RANDOM FIELDS ON THE PLANE[1]

By Ethan B. Anderes and Michael L. Stein

*University of California, Berkeley and University of Chicago*

This paper presents a new approach to the estimation of the deformation of an isotropic Gaussian random field on $\mathbb{R}^2$ based on dense observations of a single realization of the deformed random field. Under this framework we investigate the identification and estimation of deformations. We then present a complete methodological package—from model assumptions to algorithmic recovery of the deformation—for the class of nonstationary processes obtained by deforming isotropic Gaussian random fields.

**1. Introduction.** Random fields obtained by deforming the coordinates of an isotropic random field form a rich class of nonstationary processes. Such random fields have the form $Y(\mathbf{x}) = Z(f^{-1}(\mathbf{x}))$, where $Z$ is a random field on $\mathbb{R}^d$ and $f$ is a deformation, or a smooth bijection of $\mathbb{R}^d$. In one dimension, for example, the deformation $f$ is used in speech recognition to model local compression or expansion of time in different utterances of a spoken word (see [17]). Working with deformations in more than one dimension, however, has been a challenge. One dimensional deformations behave locally as a change of scale. In more than one dimension, however, deformations can rotate as well as scale local coordinates. These added dynamics, in addition to complicated requirements for invertibility, have restricted the pragmatic use of deformations when modeling processes in more than one dimension. In this paper we establish methodology for working with and estimating

Received October 2005; revised November 2006.

[1]Supported by U.S. EPA—Science To Achieve Results (STAR) Program, Grant R-82940201. Although the research described in this article has been funded wholly by the United States Environmental Protection Agency, it has not been subjected to the Agency's required peer and policy review and therefore does not necessarily reflect the views of the Agency, and no official endorsement should be inferred.

deformations in $\mathbb{R}^2$ from a single realization of a deformed isotropic Gaussian random field.

The use of deformations to model nonstationary processes was first introduced to the spatial statistics literature by Sampson and Guttorp [18]. Their work, as well as that of subsequent authors (see, e.g., [7, 11, 16, 19]) deals mostly with data from sparse observation locations with independent replicates of the random field. Two recent papers by Clerc and Mallat [4, 5] consider a similar deformation model, but under a different observation scenario: a densely observed single realization. They use families of localized functions to estimate local properties of the deformation, and under the assumption of reflective shape recovery, these local properties are related to estimates of the shape of the reflective surface. Most of their results hold in one dimension, however, and it is not clear that in two dimensions their estimates work for arbitrarily smooth isotropic processes under general deformations.

Guyon and Perrin [10] tackle the problem of developing consistent estimates of deformations in two dimensions. They succeed in proving results within a class of deformations when observing random fields that are stationary but not isotropic. In this paper we examine the consequence of adding the assumption of isotropy to the pre-deformed random fields. The isotropic assumption complicates the estimation of the original orientation of the local coordinates by introducing a rotational invariance to the random field. In Section 3 we notice the local behavior of $C^1$ diffeomorphisms can be approximated by a composition of rotations, linear coordinate stretchings and translations. Estimating the initial local rotation—or the original orientation of the local coordinates—using only observations in a local neighborhood becomes difficult. Our approach is to estimate two local parameters of the deformation, a dilatation and scale, which are invariant under initial local rotations. These parameters can be estimated from the local behavior of the deformed process and are still sufficient to uniquely characterize the deformation, up to global rigid motions.

The structure of this paper is as follows. In Section 2 we present our modeling assumptions on the random fields to be deformed. We introduce a flexible semi-parametric class of random fields under which we can perform approximate likelihood estimation of the local parameters that characterize the deformation. In Section 3 we make explicit our assumptions on the class of deformations and study in more detail the consequences of isotropy when estimating deformations. Sections 4 and 5 present the main contribution of this paper. We outline our methodology and present an algorithm which allows fast construction of the deformation. Finally, we discuss some simulations and future work.

**2. Deformation model.** Our deformation model for nonstationarity is processes of the form $Z \circ f^{-1}$, where $Z$ is a constant mean isotropic Gaussian process on $\mathbb{R}^2$ and $f:\mathbb{R}^2 \to \mathbb{R}^2$ is an orientation preserving $C^1$ diffeomorphism. We also add the assumption $f$ has bounded local distortion, which implies our maps are quasiconformal (see Appendix 6 for the definition of quasiconformal and $C^1$ diffeomorphism). By bounded local distortion, we mean that the ratio of $\limsup_{\mathbf{x} \to \mathbf{x}_0} |f(\mathbf{x}) - f(\mathbf{x}_0)|/|\mathbf{x} - \mathbf{x}_0|$ to $\liminf_{\mathbf{x} \to \mathbf{x}_0} |f(\mathbf{x}) - f(\mathbf{x}_0)|/|\mathbf{x} - \mathbf{x}_0|$ be uniformly bounded for all $\mathbf{x}_0 \in \mathbb{R}^2$.

One of the advantages of this model is that estimates of $f^{-1}$ can be used to transform the coordinates of $Z \circ f^{-1}$, returning a nonstationary process to isotropy for which one can use existing statistical techniques. This can be done by mapping the observation locations $\mathbf{x}_i$ to $f^{-1}(\mathbf{x}_i)$, which transforms the graph of $(\mathbf{x}_i, Z \circ f^{-1}(\mathbf{x}_i))$ to the graph of the original isotropic process $Z$ at the locations $f^{-1}(\mathbf{x}_i)$.

We want the second order behavior of $Z$ to be as general as possible while still allowing approximate likelihood techniques for estimation. We do this by introducing a regularization parameter $\alpha > 0$ which restricts the behavior of the autocovariance at the origin, and therefore controls the mean square smoothness of the isotropic processes $Z$. To simplify notation, let $p_\alpha$ denote the nonnegative integer

$$(1) \qquad p_\alpha = \begin{cases} \lfloor \alpha/2 \rfloor, & \text{if } \alpha/2 \notin \mathbb{Z}, \\ \alpha/2 - 1, & \text{if } \alpha/2 \in \mathbb{Z}. \end{cases}$$

Let $K$ be the autocovariance function for $Z$ so that $K(|\mathbf{t}-\mathbf{s}|) = \text{cov}\{Z(\mathbf{t}), Z(\mathbf{s})\}$ for $\mathbf{s}, \mathbf{t} \in \mathbb{R}^2$. We suppose there exists a constant $c > 0$ and an $\alpha > 0$ so that $K(|t|)$ has $2p_\alpha$ continuous derivatives and

$$(2) \qquad K(|t|) - \sum_{k=0}^{p_\alpha} \frac{K^{(2k)}(0)}{(2k)!} t^{2k} \sim c\, G_{\alpha(t)} \qquad \text{as } |t| \to 0,$$

where $G_\alpha$ is defined by

$$G_\alpha(t) = \begin{cases} (-1)^{1+\lfloor \alpha/2 \rfloor} |t|^\alpha, & \text{for } \alpha/2 \notin \mathbb{Z}, \\ (-1)^{1+\alpha/2} |t|^\alpha \log |t|, & \text{for } \alpha/2 \in \mathbb{Z}. \end{cases}$$

This class of autocovariances encompasses a broad range of processes including the Matérn model (see [21]) and the so-called exponential family $\exp(-c|\mathbf{t}-\mathbf{s}|^\gamma)$ for $\gamma \in (0,2)$. The parameter $\alpha$ controls the mean square differentiability of $Z$ so that $Z$ is $n$ times mean square differentiable if and only if $n < \alpha/2$.

One advantage of this class of processes is that we can perform restricted maximum likelihood estimation under the distributional approximation of $Z$ by an intrinsic random function of order $\lfloor \alpha/2 \rfloor$ with generalized covariance function $cG_\alpha$ (see [21]). In particular, suppose $Z$ is a process with autocovariance function $K$ that satisfies (2) and suppose we observe $\mathbf{Z} := (Z(\mathbf{x}_1), \ldots, Z(\mathbf{x}_n))$. We can find a matrix $\mathbf{L}$ with $n$ columns so that the covariance structure of $\mathbf{LZ}$ is easily approximated by a positive definite matrix that depends only on the function $cG_\alpha$. To find such a matrix $\mathbf{L}$, we first need some notation: for $\mathbf{x} = (x, y) \in \mathbb{R}^2$ and $\mathbf{r} = (r_1, r_2) \in \mathbb{N}^2$, let $\mathbf{x}^\mathbf{r}$ denote the monomial $x^{r_1} y^{r_2}$. Now let $\mathbf{L}$ be a matrix for which each row $(u_1, \ldots, u_n)$ satisfies $\sum_{i=1}^n u_i \mathbf{x}_i^\mathbf{r} = 0$ for all $\mathbf{r} = (r_1, r_2) \in \mathbb{N}$ such that $r_1 + r_2 \leq \lfloor \alpha/2 \rfloor$. These rows can be easily computed by finding linearly independent vectors orthogonal to the space spanned by $\{(\mathbf{x}_1^\mathbf{r}, \ldots, \mathbf{x}_n^\mathbf{r}) : r_1 + r_2 \leq \lfloor \alpha/2 \rfloor\}$. Now let $(u_1, \ldots, u_n)$ and $(v_1, \ldots, v_n)$ be two rows of $\mathbf{L}$ and let $k \leq \lfloor \alpha/2 \rfloor$, then

$$\sum_{i=1}^n \sum_{j=1}^n u_i v_j |\mathbf{x}_i - \mathbf{x}_j|^{2k} = 0$$

so that

$$\text{(3)} \quad \text{cov}\left\{\sum_{i=1}^n u_i Z(\mathbf{x}_i), \sum_{j=1}^n v_j Z(\mathbf{x}_j)\right\}$$
$$= c \sum_{i,j=1}^n u_i v_j G_\alpha(|\mathbf{x}_i - \mathbf{x}_j|) + \sum_{i,j=1}^n u_i v_j R(|\mathbf{x}_i - \mathbf{x}_j|),$$

where $R$ denotes the error term in the asymptotic approximation (2) so that $R(|t|) = o(G_\alpha(|t|))$ as $|t| \to 0$. Let $(u_1^\ell, \ldots, u_n^\ell)$ denote the $\ell$th row of $\mathbf{L}$. Since the functions $G_\alpha$ are conditionally positive definite of order $\lfloor \alpha/2 \rfloor$, the matrix $(c \sum_{i,j=1}^n u_i^p u_j^q G_\alpha(|\mathbf{x}_i - \mathbf{x}_j|))_{p,q}$ is positive definite (see [21]). Therefore, by ignoring the last term in (3), we can approximate the covariance structure of $\mathbf{LZ}$, under the semiparametric assumption (2) by the covariance matrix $(c \sum_{i,j=1}^n u_i^p u_j^q G_\alpha(|\mathbf{x}_i - \mathbf{x}_j|))_{p,q}$.

We point out a useful fact that will be used in the following methodology: if one supposes the deformation $f$ is smoother than the sample paths of the process $Z$, then $Z$ and $Z \circ f^{-1}$ share the same local mean squared smoothness. This is to our advantage in that $\alpha$ is invariant under deformations of the process, which should allow estimation of $\alpha$ directly from $Z \circ f^{-1}$ without knowledge of $f$. For example, suppose the autocovariance of the process $Z$ has the form (2) with $\alpha < 2$ and suppose $\mathbf{h} \in \mathbb{R}^2$ such that $\mathbf{h} \neq (0, 0)$. We show that there exists a constant $0 < c_1 < \infty$ such that

$$\text{(4)} \quad \mathbb{E}\{Z(\mathbf{x} + \epsilon\mathbf{h}) - Z(\mathbf{x})\}^2 \sim c_1 \mathbb{E}\{Z \circ f^{-1}(\mathbf{x} + \epsilon\mathbf{h}) - Z \circ f^{-1}(\mathbf{x})\}^2$$

as $\epsilon \to 0$. In other words, the processes $Z$ and $Z \circ f^{-1}$ both have the same power decay of the variogram at the origin (see [6]), which can be estimated from one realization of $Z \circ f^{-1}$.

To see why (4) is true, first notice that

$$\text{(5)} \quad \mathbb{E}\{Z(\mathbf{x} + \epsilon\mathbf{h}) - Z(\mathbf{x})\}^2 = 2K(0) - 2K(|\epsilon\mathbf{h}|) \sim 2c|\mathbf{h}|^\alpha \epsilon^\alpha$$

as $\epsilon \to 0$ (for the rest of this section we write $f_1 \sim f_2$ to denote that $f_1/f_2 \to 1$ as $\epsilon \to 0$ with $\mathbf{x}$ and $\mathbf{h}$ fixed). Similarly,

$$\mathbb{E}\{Z \circ f^{-1}(\mathbf{x}+\epsilon\mathbf{h}) - Z \circ f^{-1}(\mathbf{x})\}^2 \sim 2c|f^{-1}(\mathbf{x}+\epsilon\mathbf{h}) - f^{-1}(\mathbf{x})|^\alpha$$

by direct application of (2) since $|f^{-1}(\mathbf{x}+\epsilon\mathbf{h}) - f^{-1}(\mathbf{x})| \to 0$ as $\epsilon \to 0$. Now by the differentiability assumption on $f^{-1}$ we also have $|f^{-1}(\mathbf{x}+\epsilon\mathbf{h}) - f^{-1}(\mathbf{x})| \sim c_2 \epsilon$ for some $0 < c_2 < \infty$. This gives

$$\mathbb{E}\{Z \circ f^{-1}(\mathbf{x}+\epsilon\mathbf{h}) - Z \circ f^{-1}(\mathbf{x})\}^2 \sim 2c\,c_2^\alpha\,\epsilon^\alpha,$$

which, in conjunction with (5), proves (4).

The following methodology requires that the process $Z$ be smoother than the deformation $f$. However, this assumption is made only to exclude difficulties when estimating $\alpha$. In addition, we anticipate only using this method in practice for smooth deformations where the process $Z$ is rarely very smooth. Therefore, we consider this assumption a minor technicality to the methods described below.

REMARK 1. In what follows we will be using approximate likelihood methods that will depend on the data only through sufficiently high order increments of the processes. Therefore, all our methods can be extended to intrinsic random functions (see [3]) such as fractional Brownian surfaces, for example.

REMARK 2. The dependence of $K$ on the parameter $c$ presents an identifiability problem if one assumes the constant $c$ is unknown. This is because $c$ is confounded with changes of scale (which are included in our deformation class). For example, if the autocovariance for $Z$ satisfies (2) with $c=1$, then the autocovariance for the deformed processes $Z(2\mathbf{x})$ satisfies (2) for $c=2^\alpha$. In what follows we choose to fix the constant $c=1$ and assume it is known. This allows us to estimate the correct scale of the deformation.

**3. $C^1$ diffeomorphisms with bounded distortion.** Since the original process, $Z$, is assumed to be isotropic, the most one can hope for is the identification of the deformation up to rigid motions of the plane. We will show in this section that this is indeed possible. We begin by studying the local behavior of the map $f$ in terms of local affine transformations—whose existence is guaranteed by the assumption that $f$ is a $C^1$ diffeomorphism. These local affine transformations can be decomposed using the singular value decomposition, which makes explicit how to obtain the local coordinates of the map by a composition of rotations and stretchings. Since the isotropic assumption on the original random field complicates the identification of the initial local rotation, we estimate the remaining parameters of the local coordinates. These remaining parameters are characterized by an ellipse. The theory of quasiconformal maps, in part, studies maps through these ellipses, and provides us with the theoretical foundation we need. A brief overview of quasiconformal maps is included in Appendix 6.

Since $f$ is a diffeomorphism, its local linear behavior is characterized by

$$f(\mathbf{x} + \mathbf{h}) = f(\mathbf{x}) + J_f \mathbf{h} + o(|\mathbf{h}|), \qquad \text{as } |\mathbf{h}| \to 0, \tag{6}$$

where $J_f$ is the Jacobian of $f$ and, writing $f(x, y) = (u(x, y), v(x, y))$, is defined as

$$J_f := \begin{pmatrix} u_x & u_y \\ v_x & v_y \end{pmatrix}. \tag{7}$$

Since $f$ is orientation-preserving, $\det J_f > 0$, and since $f$ is $C^1$, it has continuous partial derivatives $u_x$, $u_y$, $v_x$ and $v_y$. We can further decompose the map $f$ by using the singular value decomposition to represent the linear map $J_f$ as a sequence of rotations and stretchings. In particular, $J_f = U\Lambda V^t$, where $U$, $V$ are orthogonal matrices and $\Lambda$ is a diagonal matrix with diagonal elements $\lambda_1 \geq \lambda_2 > 0$ (see [9]). Since $\det J_f > 0$, we can take $U, V$ to be rotation matrices, which gives

$$J_f = \begin{pmatrix} \cos\theta & -\sin\theta \\ \sin\theta & \cos\theta \end{pmatrix} \begin{pmatrix} \lambda_1 & 0 \\ 0 & \lambda_2 \end{pmatrix} \begin{pmatrix} \cos\psi & -\sin\psi \\ \sin\psi & \cos\psi \end{pmatrix}, \tag{8}$$

where $\lambda_1 \geq \lambda_2 > 0$. This decomposition describes the local behavior of $f$ by an initial rotation, a stretching and a final rotation. This sheds light on the information $Z \circ f^{-1}$ provides about the map $f$. Since $Z$ is rotationally invariant, the initial rotation by an angle of $\psi$ will be difficult, if not impossible, to estimate locally. However, our deformations map infinitesimal circles to infinitesimal ellipses and these ellipses are invariant under the initial local rotation of the map $f$, so there is some prospect of estimating them from one sample path of the process $Z \circ f^{-1}$. The inclinations of these ellipses are given by $\theta$ and their eccentricities by $\lambda_1/\lambda_2 \geq 1$.

These locally defined ellipses are particularly useful in the study of deformations. First, notice that the eccentricity provides a natural measure of local distortion that equals our original notion of distortion as the ratio of the upper and lower limits of $|f(\mathbf{x}) - f(\mathbf{x}_0)|/|\mathbf{x} - \mathbf{x}_0|$ as $\mathbf{x}$ approaches $\mathbf{x}_0$ from different directions. This notion of distortion will become useful when studying the smoothing problem for local ellipse estimation. Second, there is existing literature that characterizes maps with a given ellipse field. For example, the inclinations and eccentricities of these ellipses are sufficient to recover $f^{-1}$ up to postcomposition of conformal maps (see Appendix 6).

We parameterize these ellipses by a pair $\{(\mu, \phi) \in \mathbb{C} \times \mathbb{R}^+ : |\mu| < 1\}$, where $\mu$ is referred to as the complex dilatation and $\phi$ as the scale. Set $\mu$ to be the complex number

$$\mu = -\frac{K-1}{K+1} e^{i2\theta}, \tag{9}$$

where $K$ denotes the eccentricity of the ellipse and $\theta$ denotes inclination. The scale parameter is defined by $\phi = \sqrt{\lambda_1 \lambda_2}$. In the mathematical literature $\mu$ is called the complex dilatation of the map $f^{-1}$. The complex dilatation of a map characterizes the infinitesimal ellipse that gets mapped to a infinitesimal circle. The mapping theorem from quasiconformal theory tells us that specifying $(\mu, \phi)$ on a Jordan region $\Omega$ is sufficient to characterize the map $f^{-1}$ on $\Omega$ up to rigid motions of $\Omega$ (see Appendix 6).

**4. Methodology.** This section outlines our procedure for the estimation of the deformation $f^{-1}$ and the fractional index of the process $\alpha$ in the presence of a densely observed sample path of a deformed Gaussian process. The estimation procedure will be done in three stages. First we estimate the fractional index of the process, taking advantage of the invariance of the fractal properties of the sample paths under sufficiently smooth deformations. Then, using the estimated fractional index $\hat{\alpha}$, we use a local likelihood approach to estimate the local complex dilatation and scale of the deformation, $(\mu, \phi)$, independently over each neighborhood. Finally, we smooth and interpolate the estimates $(\hat{\mu}, \hat{\phi})$ across neighborhoods.

The two parameters $(\mu, \phi)$ are enough to characterize the deformation $f^{-1}$ on the observation region up to a rotation and translation of the original coordinates. Therefore, on estimating $(\mu, \phi)$, one has estimated all the parameters in the model. To make this more useful, however, we need to recover $f^{-1}$ from $(\mu, \phi)$. Sections 4.3.2 and 5 detail how one can efficiently reconstruct the deformation given $(\mu, \phi)$.

For both the estimation of $\alpha$ and $(\mu, \phi)$, we use likelihood methods. Since we are supposing we have observed the process densely, any full likelihood method will be computationally impractical. Therefore, we partition the observation locations into neighborhoods and assume independence of the process across partitions. To be clear, we are not changing our model assumptions, just devising approximate likelihood techniques. Likelihood methods present two advantages for an initial study of estimation. First, there are no requirements on the configuration of observation locations. Second, the estimates are easily constructible and our experience has been that they are highly efficient.

It is advantageous to switch to complex notation where a point $(x, y)$ is represented by $z = x + iy$ and consider our deformations as mapping regions of the complex plane. To set notation, we let $\mathbf{z} := (z_1, \ldots, z_n)'$ be the vector of observation locations in $\mathbb{C}$, and $\mathbf{Y} := (Y_1, \ldots, Y_n)'$ denote the corresponding observations of the process $Y := Z \circ f^{-1}$. In what follows we dividing the observations into local neighborhoods. We denote by $\mathcal{N}_n$ the partition of indices $\{1, \ldots, n\}$ that corresponds to the local observation neighborhoods. For an index set $\mathcal{I} = \{i_1, \ldots, i_m\} \in \mathcal{N}_n$, let $\mathbf{z}_\mathcal{I}$ denote the vector $(z_{i_1}, \ldots, z_{i_m})$, and similarly, define $\mathbf{Y}_\mathcal{I}$. Last, for a function $F(z, w)$, we let $(F(z_p, z_q))_\mathcal{I}$ denote the matrix with entries $F(z_p, z_q)$ for $p, q \in \mathcal{I}$.

The description of the methodology is described as it applies to the simulation pictured in Figure 1. It shows a simulation of a deformed process $Y := Z \circ f^{-1}$, where the isotropic autocovariance function for $Z$ has the form $K(t) = 0.5151 - |t|^{0.7} + O(|t|^2)$ as $t \to 0$; see [2] for details. The observations are on a grid of size $400 \times 400$ in $[0, 1]^2$ and the deformation that transforms the coordinates of $Y$ to isotropy is

$$f^{-1}(x + iy) := (1.2 - y)\, e^{-i\pi(1-x)/2} + i1.2. \tag{10}$$

4.1. *Estimating $\alpha$.* To estimate $\alpha$, we first construct the neighborhood structure, $\mathcal{N}_n$, of the observation locations. For this step alone, one would want the size of the neighborhoods to be as large as computationally feasible. Then $\hat\alpha$ is obtained by maximizing an approximated log-likelihood under the supposition of independence across neighborhoods and that the process is isotropic with autocovariance satisfying (2).

To approximate the log-likelihood for the observations $\mathbf{Y}_\mathcal{I}$ on a particular neighborhood $\mathcal{I} = \{i_1, \ldots, i_m\} \in \mathcal{N}_n$, we use the techniques discussed in Section 2 and approximate the log-likelihood for a linear transformation of the data $\mathbf{L}\mathbf{Y}_\mathcal{I}$. We want to set the rows of $\mathbf{L}$ to be linearly independent vectors orthogonal to the space spanned by $\{(z_{i_1}^\mathbf{r}, \ldots, z_{i_m}^\mathbf{r}) : r_1 + r_2 \leq \lfloor \alpha/2 \rfloor\}$, where $(x + iy)^\mathbf{r}$ denotes the real number $x^{r_1} y^{r_2}$. Unfortunately this involves knowing $\lfloor \alpha/2 \rfloor$. However, if we suppose we know a predetermined upper bound on the magnitude of $\alpha$, which we denote by $\alpha_o$, it will be sufficient to find linearly independent vectors orthogonal to space spanned by $\{(z^{\mathbf{r}}_{i_1}, \ldots, z^{\mathbf{r}}_{i_m}) : r_1 + r_2 \leq \lfloor \alpha_o/2 \rfloor\}$. Now, using (3), the random vector $\tilde{\mathbf{Y}}_\mathcal{I} := \mathbf{L}\mathbf{Y}_\mathcal{I}$ is multivariate Gaussian with covariance matrix

$$E\tilde{\mathbf{Y}}_\mathcal{I}\tilde{\mathbf{Y}}'_\mathcal{I} = \mathbf{L}(K(|z_p - z_q|))_\mathcal{I}\mathbf{L}'$$
$$= \mathbf{L}(G_\alpha(|z_p - z_q|))_\mathcal{I}\mathbf{L}' + \mathbf{L}(R(|z_p - z_q|))_\mathcal{I}\mathbf{L}'.$$

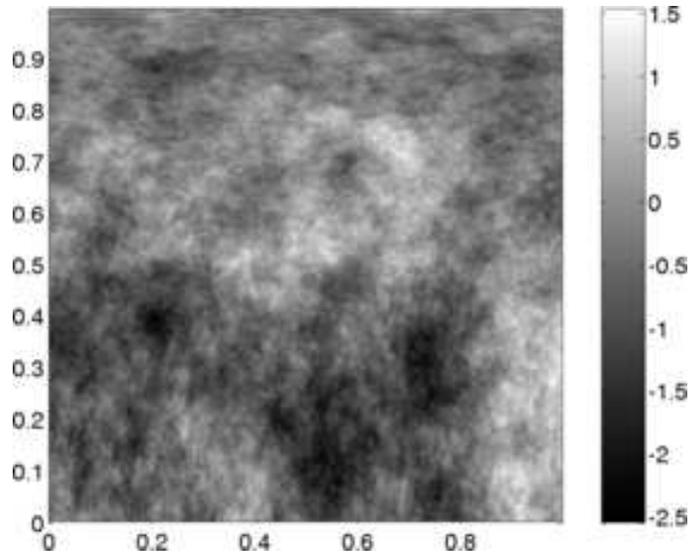

FIG. 1. *Deformed process.*

By ignoring the remainder term $\mathbf{L}(R(|z_p - z_q|))_\mathcal{I}\mathbf{L}'$ and setting $\Sigma_{\mathcal{I},\alpha} := \mathbf{L}(G_\alpha(|z_p - z_q|))_\mathcal{I}\mathbf{L}'$, we get the following approximate log likelihood over neighborhood $\mathcal{I}$:

(11) $$\ell(\alpha|\tilde{\mathbf{Y}}_\mathcal{I}, \mathbf{z}_\mathcal{I}) = -\tfrac{1}{2}\log|\Sigma_{\mathcal{I},\alpha}| - \tfrac{1}{2}\tilde{\mathbf{Y}}'_\mathcal{I}\Sigma^{-1}_{\mathcal{I},\alpha}\tilde{\mathbf{Y}}_\mathcal{I} + C,$$

where $C$ is a constant additive term not depending on $\alpha$. We suppose independence across blocks so the log-likelihood for the full observation vector $\mathbf{Y}$ is computed by summing the log-likelihoods (11) over neighborhoods $\mathcal{I} \in \mathcal{N}$. Now we estimate $\alpha$ by

$$\hat{\alpha} := \underset{\alpha \leq \alpha_o}{\arg\max} \sum_{\mathcal{I} \in \mathcal{N}_n} \ell(\alpha|\tilde{\mathbf{Y}}_\mathcal{I}, \mathbf{z}_\mathcal{I}).$$

For the simulation in Figure 1, the neighborhood structure, $\mathcal{N}_n$, was constructed by dividing the $400 \times 400$ grid into $40 \times 40$ blocks of size $10 \times 10$. The maximum approximate likelihood estimate is $\hat{\alpha} = 0.696$ (the true $\alpha$ is 0.7). We do not address the accuracy of the estimates of $\alpha$ in this paper. However, in all our simulations we have seen very accurate estimation of $\alpha$.

4.2. *Estimating* $(\mu, \phi)$. In this section we suppose we have an estimate $\hat{\alpha}$ so that by regarding this estimate as the truth, the only unknown parameter is the quasiconformal map $f$, or equivalently, $f^{-1}$. As was discussed in Section 3, we reparameterize $f^{-1}$ by its complex dilatation and scale field $(\mu, \phi)$. Fixing $\hat{\alpha}$, we use approximate maximum likelihood estimation on $(\mu, \phi)$ locally on each neighborhood to get a spatial map $(\hat{\mu}, \hat{\phi})$.

Given a particular neighborhood $\mathcal{I} \in \mathcal{N}_n$, we suppose the process behaves like a geometric anisotropic process $Z(J_{f^{-1}}\mathbf{x})$ with smoothness parameter $\hat{\alpha}$. As was discussed in Section 3, only three parameters of $J_{f^{-1}}$ are estimable, and they can be reparameterized as a complex dilatation and scale $(\mu_i, \phi_i)$. In complex notation the geometric anisotropic process can be written as $Z(A_i(z + \mu_i\overline{z}))$, where $|A_i|^2 = \phi_i^2/(1-|\mu_i|^2)$. Here $i$ indexes the neighborhoods $\mathcal{N}_n$ and $(\mu_i, \phi_i)$ denotes these constants for each neighborhood. As in Section 4.1 consider the vector $\tilde{\mathbf{Y}}_\mathcal{I} := \mathbf{L}\mathbf{Y}_\mathcal{I}$, where the rows of $\mathbf{L}$ are linearly independent vectors orthogonal to $\{(z^{\mathbf{r}}_{i_1}, \ldots, z^{\mathbf{r}}_{i_m}) : r_1 + r_2 \leq \lfloor \hat{\alpha}/2 \rfloor\}$. Letting $\theta_i$ denote $(\mu_i, \phi_i)$, we then estimate $\theta_i$ by maximizing the log-likelihood $\ell(\theta_i|\tilde{\mathbf{Y}}_I, \mathbf{z}_\mathcal{I}) = -\log|\Sigma_{\mathcal{I},\theta_i}|/2 - \tilde{\mathbf{Y}}'_\mathcal{I}\Sigma^{-1}_{\mathcal{I},\theta_i}\tilde{\mathbf{Y}}_\mathcal{I}/2 + C$, where $C$ is a constant

term not depending on $\theta_i$. Similar to what was done in Section 4.1, the covariance matrix $\Sigma_{\mathcal{I},\theta_i}$ is computed by ignoring the remainder term so that $\Sigma_{\mathcal{I},\theta_i} = \mathbf{L}(G_{\hat{\alpha},\theta_i}(z_p - z_q))_{\mathcal{I}}\mathbf{L}'$, where $G_{\hat{\alpha},\theta_i}(z) := G_{\hat{\alpha}}(|A_i|(z - \mu_i \overline{z}))$. Notice that the parameter range for $\mu$ is the complex unit disk $\{\mu \in \mathbb{C} : |\mu| < 1\}$ and $\phi > 0$.

For the simulation in Figure 1, we again use the neighborhood structure of $40 \times 40$ blocks of size $10 \times 10$ to estimate $(\mu, \phi)$. The left-hand plot in Figure 2 graphs the local estimates $\mu_i$ as vectors.

4.3. *Smoothing and interpolating* $(\hat{\mu}, \hat{\phi})$. The parameters $(\hat{\mu}, \hat{\phi})$ are locally estimated values of the functions $\mu$ and $\phi$. Since the local likelihood estimation was done independently on each neighborhood, there was no smoothness constraint incorporated into the estimation of $\mu$ and $\phi$. To incorporate smoothness conditions in these estimates, we choose to smooth $(\hat{\mu}, \hat{\phi})$ after local likelihood estimation has been completed. As we will see, smoothing $\hat{\mu}$ and $\hat{\phi}$ should, and will, be done in completely different ways.

4.3.1. *The complex dilatation field* $\hat{\mu}$. The most naïve approach to smoothing the dilatation field $\hat{\mu}$ would be local averaging. However, it is worthwhile to investigate the smoothing problem under more generality. In particular, the dilatation field is a parameterization of the ellipse, which is the pre-image of a local infinitesimal circle under the deformation $f^{-1}$. Just averaging the elements of a particular parameterization, however, does not provide a consistent notion of smoothing under reparameterization. In the following we develop a metric on the complex dilatation parameterization, based on a geometrical motivation, and use this metric to smooth $\hat{\mu}$.

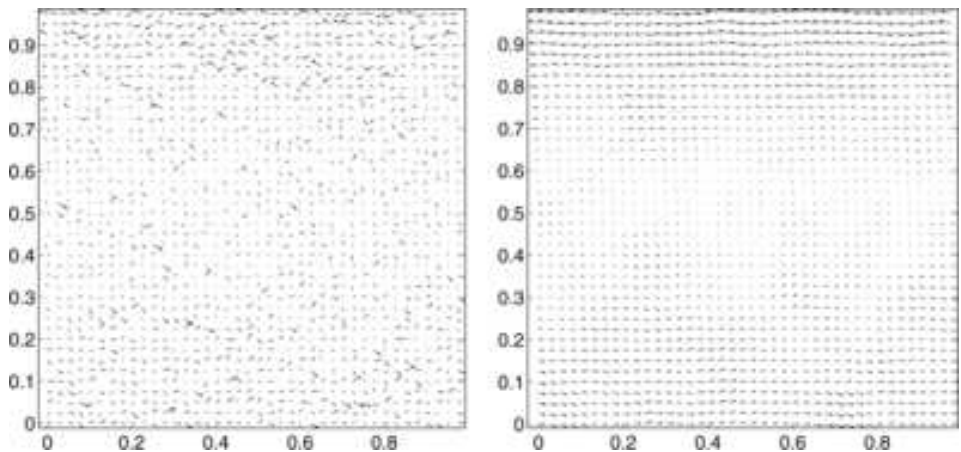

FIG. 2. *Left*: *The estimated complex dilatation* $\hat{\mu}$. *Right*: *The result of smoothing the estimated dilatation* $\hat{\mu}$.

The locally defined ellipses represent infinitesimal linear deformations, the magnitudes of which are naturally measured by their eccentricities, denoted by $K$. To develop a metric on ellipses, consider two matrix representations $A$ and $B$ of infinitesimal linear transformations. We claim that the magnitude of the distortion of the composed linear map $AB^{-1}$ is a natural measure of distance between two ellipses. By (9), this is easily computed by

$$K_{AB^{-1}} = \frac{1 + |\mu_{AB^{-1}}|}{1 - |\mu_{AB^{-1}}|}, \qquad |\mu_{AB^{-1}}| = \left|\frac{\mu_A - \mu_B}{1 - \mu_A \overline{\mu}_B}\right|.$$

Since $K_{AB^{-1}}$ is at least 1, we take logarithms to define the distance measure $d$ on the complex dilatations:

$$d(\mu_A, \mu_B) := \frac{1}{2} \log K_{AB^{-1}} = \frac{1}{2} \log\left(\frac{1 + |\mu_{AB^{-1}}|}{1 - |\mu_{AB^{-1}}|}\right).$$

It interesting to note that $d$ is indeed a metric and it equals the hyperbolic metric on the unit disk of $\mathbb{C}$ (see [12], e.g.). Much is known about this metric, in particular, the geodesic connecting two points is easily expressible as the arc on a circle that is orthogonal to the unit circle and joins the two points. One way to extend these ideas to define the convex combination, $\mu^*$, of $n$ complex dilatations with positive weights $w_1, \ldots, w_n$ could be

$$(12) \qquad \mu^* = \arg \min_{|\mu| < 1} \{w_1 \, d^p(\mu, \mu_1) + \cdots + w_n \, d^p(\mu, \mu_n)\}$$

for some power $p > 0$. Now by using weights that spatially vary, one can smooth the estimated dilatation field $\hat{\mu}$. In addition, convex combinations allow one to define bilinear or other interpolation approaches to produce maps of $\hat{\mu}$ throughout the observation region $\Omega$.

Unfortunately, we do not know of any closed form for (12), but a computational approximation to the minimum is attainable. The right-hand plot of Figure 2 shows the result of the approximated minimization of (12) for $p = 2$ with uniform weights on a sliding window of size $4 \times 4$ in the interior, and the rectangular region that overlaps the sliding $4 \times 4$ window near the boundary.

4.3.2. *The scale field $\hat{\phi}$.* Given a quasiconformal map $f: \Omega \to \Omega'$ with complex dilatation $\mu$, all other quasiconformal maps with dilatation $\mu$ are obtained by postcomposing $f$ with a conformal map on $\Omega'$. Therefore, once we smooth and interpolate $\hat{\mu}$, we need to use the scale information $\phi$ to identify the deformation within the class of quasiconformal maps specified by $\hat{\mu}$. Section 5 outlines an algorithm for constructing a representative quasiconformal map with prescribed dilatation. We use this map to transform the observation coordinates, then use the scale information to postcompose this map conformally to estimate the true deformation $f^{-1}$. Figure 3 shows an example of two maps with the same dilatation but different scale fields. The right-hand plot is the map $f^{-1}$ given by equation (10) and the left-hand plot is a representative map, constructed from our algorithm, with the same dilatation $\mu_{f^{-1}}$. The two maps can be equated by postcomposing one with a conformal map.

Consider deforming the observation space $\Omega$ by a representative quasi-conformal map $\check{f}: \Omega \to \Omega'$ with dilatation $\mu^*$ (see Section 5), where $\mu^*$ is the dilatation obtained by smoothing $\hat{\mu}$. Remember that $\mu^*$ is an estimate of $\mu_{f^{-1}}$ so that $\check{f}$ is an estimate—before postcomposition with a conformal map—of $f^{-1}$. So by transforming the coordinates of the observed process $Z \circ f^{-1}$ by $\check{f}$, the resulting process takes the form $Z \circ g^{-1}$, where $g = \check{f} \circ f$ is a conformal map with local scale $\phi_g = |g'| = \phi_f(\phi_{\check{f}} \circ f)$. Therefore, we want to find a conformal $h$ that transforms $g$ to the identity, that is, $h = g^{-1}$. Once we find such an $h$, the composition $h \circ \check{f}$ will estimate $f^{-1}$ up to rigid motions of $\mathbb{C}$. In particular, we want $|(h \circ g)'| = |h' \circ g||g'| = 1$, so that

$$(13) \qquad |h'| = \frac{1}{\phi_g \circ g^{-1}} = \frac{1}{\phi_f(\phi_{\check{f}} \circ f)} \circ g^{-1}$$

on $\Omega'$. Notice that our estimate $\hat{\phi}$ is a pointwise estimate of $\phi_{f^{-1}} = 1/\phi_f \circ f^{-1}$ at spatial locations $\mathbf{z} = (z_1, \ldots, z_n)$. Therefore, we want conformal $h$ that satisfies

$$(14) \qquad |h'(w_j)| \approx \frac{\hat{\phi}_j}{\phi_{\check{f}} \circ \check{f}^{-1}(w_j)},$$

where $(w_1, \ldots, w_n) := (\check{f}(z_1), \ldots, \check{f}(z_n))$ are points in $\Omega'$. Section 5 shows how one constructs $\check{f}$ and computes $\phi_{\check{f}} \circ \check{f}^{-1}$ efficiently.

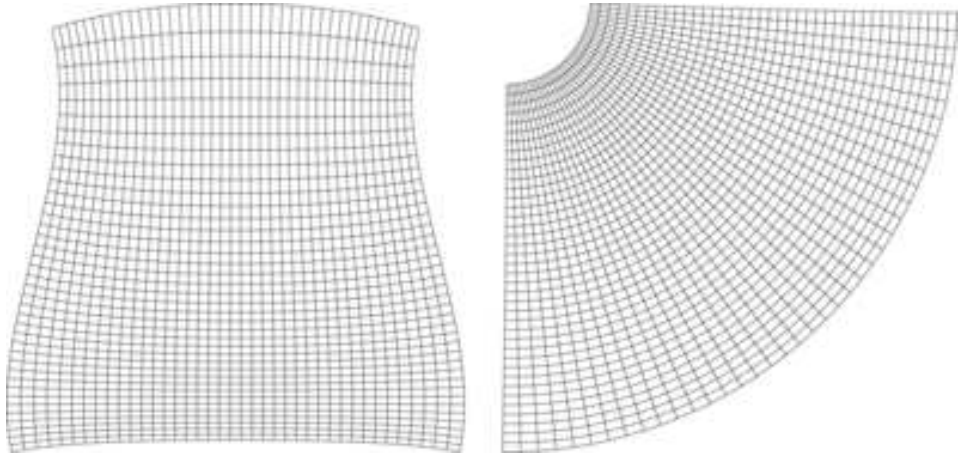

FIG. 3. *Two maps which have the same local dilatation but different local scale.*

Notice that $h'$ is analytic on $\Omega'$ and nonzero by invertibility. Therefore, $\log h'$ is analytic on $\Omega'$, with real part $\log|h'|$. Since the real part of an analytic function is harmonic, by (14), we have that $\log \hat{\phi}_j - \log \phi_{\check{f}} \circ \check{f}^{-1}(w_j)$ should approximate a sampled harmonic function at locations $w_j$. By rescaling, we suppose, without loss of generality, that $\Omega$ can be embedded into the unit disk $\mathbb{D}$ of the complex plane. Now, any harmonic function on the unit disk is harmonic when restricted to the region $\Omega \subset \mathbb{D}$. We then use the Fourier series to represent the boundary values for harmonic functions on $\mathbb{D}$. Since the boundary values of a harmonic function uniquely determine its values on the interior, this gives a series representation of harmonic functions on $\Omega$. Unfortunately there are harmonic functions on $\Omega$ that cannot be extended harmonically to the whole disk $\mathbb{D}$. One way to avoid this difficulty is to conformally map $\Omega$ to the disk $\mathbb{D}$ instead of embedding it. This is guaranteed to work, but is more computationally expensive. For this reason, our examples used the embedding technique.

We consider the class of harmonic functions on the unit disk represented by $g(w) = \sum_{n=0}^{N} A_n w^n$ for complex coefficients $A_n$ and some finite $N$. In polar coordinates $(r, \theta)$, where $w = re^{i\theta}$, this becomes $g(r, \theta) = \sum_{n=0}^{N} A_n r^n e^{in\theta}$. We use this decomposition for the harmonic function $\log h'$ on $\Omega'$. Since the real part of $\log h'$ is $\log|h'|$, we want to find a sequence $A_n$ such that

$$\operatorname{Re} \sum_{n=0}^{N} A_n r_j^n e^{in\theta_j} \approx \log \frac{\hat{\phi}_j}{\phi_{\check{f}} \circ \check{f}^{-1}(w_j)}$$

for $j = 1, \ldots, n$, where $w_j = r_j e^{i\theta_j}$. Now by setting $A_n = a_n + ib_n$, the problem becomes finding $a_n, b_n$ so that

$$\sum_{n=0}^{N} r_j^n (a_n \cos n\theta_j - b_n \sin n\theta_j) \approx \log \frac{\hat{\phi}_j}{\phi_{\check{f}} \circ \check{f}^{-1}(w_j)}.$$

The advantage of this representation is the linearity in the parameters $a_n$ and $b_n$. Using $\ell^2$ minimization to estimate $a_n + ib_n$, this corresponds to solving the linear problem $\min_{\mathbf{c}} \|\mathbf{Fc} - \mathbf{l}\|_2^2$, where

$$\mathbf{c} := (a_1, \ldots, a_N, b_1, \ldots, b_N),$$
$$\mathbf{l} := \left( \log \frac{\hat{\phi}_j}{\phi_{\check{f}} \circ \check{f}^{-1}(w_j)} \right)_j,$$
$$\mathbf{F} := [(r_j^n \cos n\,\theta_j)_{jn}, (-r_j^n \sin n\,\theta_j)_{jn}].$$

The right-hand plot of Figure 4 shows the resulting smoothed $\hat{\phi}$ for the simulation in Figure 1.

Once we have estimated the sequence $A_n$, we can then approximate $h'$ by $\widehat{h'}(w) = \exp\{\sum_{n=0}^{N} A_n w^n\}$, which is then used to construct $\hat{h}$. Finally, $f^{-1}$ is estimated by $\widehat{f^{-1}} := \hat{h} \circ \check{f}$. Figure 5 plots $\widehat{f^{-1}}$ (left) and the true deformation $f^{-1}$. Now, using our estimated deformation, we can transform the coordinates of $Y$ attempting to return the process to isotropy. The resulting process is shown in Figure 6.

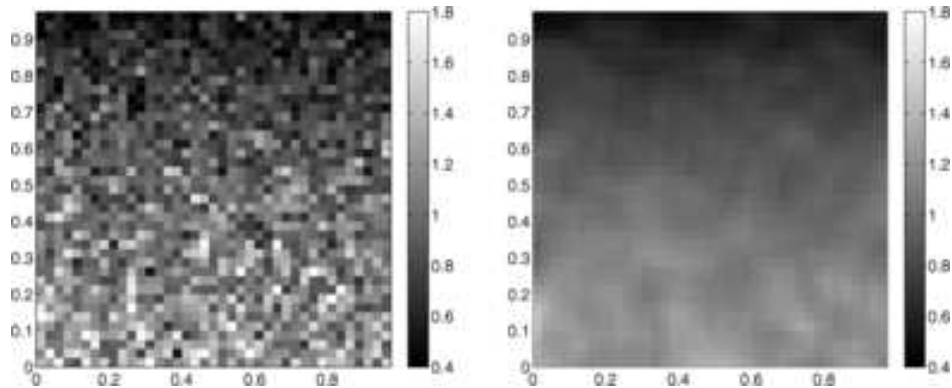

FIG. 4. *Left*: *The estimated scale field $\hat{\phi}$. Right*: *The result of smoothing the estimated scale field $\hat{\phi}$.*

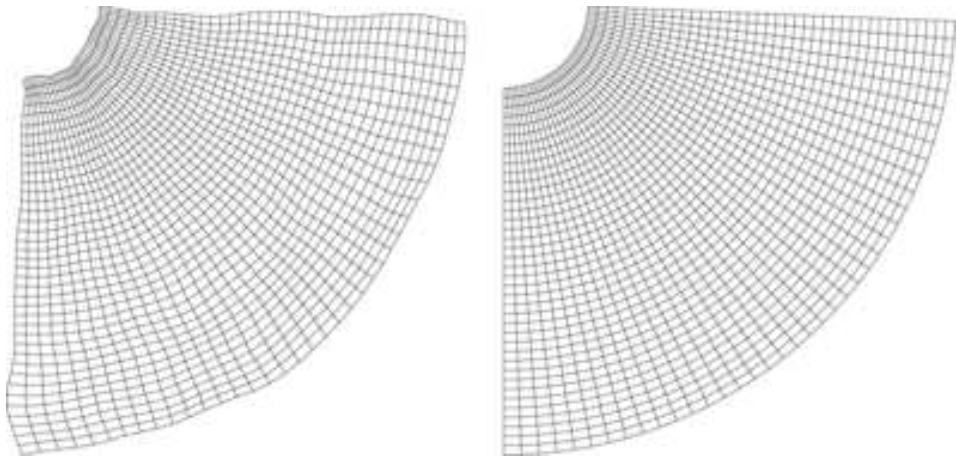

FIG. 5. *Left*: *Estimated deformation. Right*: *True deformation.*

**5. Generating a quasiconformal map from its complex dilatation.** We first need to introduce some notation and present a few facts, most of which are explained in more detail in Appendix 6. The basic approach to our algorithm is to find a vector field flow representation of a quasiconformal map. By a vector field flow representation of a quasiconformal map $f$, we mean a class of maps indexed by a time variable, $\{f_t\}_{t \in [0,1]}$, where $f_0$ is the identity, $f_1 = f$, and the functions $f_t$ depend smoothly on $t$ so that

(15) $$f_{t+\varepsilon} \circ f_t^{-1}(z) = z + \varepsilon\left(u_t(z) + iv_t(z)\right) + o(\varepsilon)$$

for some class of functions $u_t, v_t: \mathbb{C} \to \mathbb{R}$. By writing $f = u + iv$, Appendix 6 shows that the complex dilatation of $f$ can be expressed as the ratio $\partial_{\bar{z}} f / \partial_z f$, where $\partial_z f := \frac{1}{2}(u_x + v_y) + \frac{i}{2}(v_x - u_y)$ and $\partial_{\bar{z}} f := \frac{1}{2}(u_x - v_y) + \frac{i}{2}(v_x + u_y)$. The final fact that we need for the presentation of the algorithm is an equation relating the complex dilatation of the composition $g \circ f$ to the complex dilatations $\mu_g$ and $\mu_f$. Using the chain rule to compute $\partial_{\bar{z}}(g \circ f)$ and $\partial_z(g \circ f)$, we get

$$\mu_{g \circ f} := \frac{\partial_{\bar{z}}(g \circ f)}{\partial_z(g \circ f)} = \frac{\mu_f + (\overline{\partial_z f}/\partial_z f)\mu_g \circ f}{1 + (\overline{\partial_z f}/\partial_z f)\mu_g \circ f}.$$

By rearranging terms,

(16) $$\mu_g \circ f = \frac{\partial_z f}{\overline{\partial_z f}} \frac{\mu_{g \circ f} - \mu_f}{1 - \mu_{g \circ f} \overline{\mu_f}}.$$

In what follows we derive a set of differential equations for the vector fields $\{(u_t, v_t)\}_{t \in [0,1]}$ which will be numerically solved to reconstruct $f$. A continuous extension of a construction in [1], page 99, demonstrates that one can, indeed, embed a quasiconformal map $f$ into a vector field flow by specifying a class of dilatations $\{\mu_t\}_{t \in [0,1]}$ with boundary conditions $\mu_0 \equiv 0$ and $\mu_1 = \mu$ that depends smoothly on $t$. Let $\{f_t\}_{t \in [0,1]}$ be a vector field flow such that $f_t$ has dilatation $\mu_t$ with condition $\mu_0 \equiv 0$ and $\mu_1 = \mu$. Let $\{(u_t, v_t)\}_{t \in [0,1]}$ denote the vector field associated with the flow $f_t$ so that equation (15) is satisfied. By setting $g = f_{t+\varepsilon} \circ f_t^{-1}$ and $f = f_t$ in the composition formula (16), one can easily compute the dilatation of $f_{t+\varepsilon} \circ f_t^{-1}$ as

$$(17) \qquad \left( \frac{\mu_{t+\varepsilon} - \mu_t}{1 - \mu_{t+\varepsilon} \overline{\mu_t}} \frac{\partial_z f_t}{\overline{\partial_z f_t}} \right) \circ f_t^{-1} = \varepsilon \left( \frac{\partial_t \mu_t}{1 - |\mu_t|^2} \frac{\partial_z f_t}{\overline{\partial_z f_t}} \right) \circ f_t^{-1} + o(\varepsilon).$$

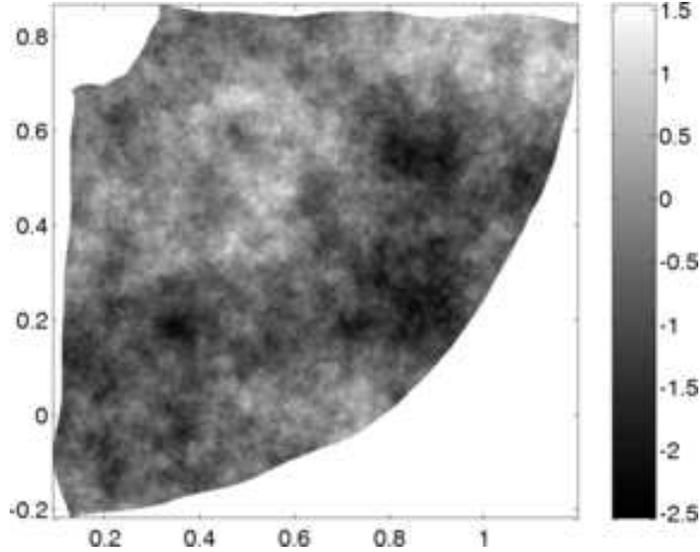

FIG. 6. *Estimated original isotropic process.*

We can also compute the complex dilatation of $z + \varepsilon(u_t + iv_t) + o(\varepsilon)$, by taking derivatives to give

$$(18) \quad \frac{\varepsilon(\partial_x u_t - \partial_y v_t) + i\varepsilon(\partial_y u_t + \partial_x v_t) + o(\varepsilon)}{2 + \varepsilon(\partial_x u_t + \partial_y v_t) + i\varepsilon(-\partial_y u_t + \partial_x v_t) + o(\varepsilon)} = \varepsilon \, \partial_{\overline{z}}(u_t + iv_t) + o(\varepsilon).$$

Then equating (17) with (18) and letting $\varepsilon \to 0$ gives

$$\partial_{\overline{z}}(u_t + iv_t) = \left( \frac{\partial_t \mu_t}{1 - |\mu_t|^2} \frac{\partial_z f_t}{\overline{\partial_z f_t}} \right) \circ f_t^{-1}.$$

In particular, the vector field $(u_t, v_t)$ satisfies

$$(19) \qquad \qquad \partial_x u_t - \partial_y v_t = 2 \operatorname{Re}(\sigma_t),$$

$$(20) \qquad \qquad \partial_y u_t + \partial_x v_t = 2 \operatorname{Im}(\sigma_t),$$

where

$$\sigma_t := \left( \frac{\partial_t \mu_t}{1 - |\mu_t|^2} \frac{\partial_z f_t}{\overline{\partial_z f_t}} \right) \circ f_t^{-1}.$$

Now, given $\{\mu_t\}_{t \in [0,1]}$, if one can solve (19) and (20) for $(u_t, v_t)$, then one can construct a time varying vector field flow realization of the map $f$. The usefulness of this representation is in its recursive nature. First notice that the initial map $f_0$ is the identity. Now supposing one has constructed $f_t$ up to some fixed time $t < 1$, one can compute $\sigma_t$ easily by deforming $\partial_t \mu_t$, $\mu_t$ and $\partial_z f_t / \overline{\partial_z f_t}$ along with $f_t$. Then by solving (19) and (20), one can approximate $f_{t+\varepsilon}$ by $f_{t+\varepsilon} = f_t + \varepsilon(u_t \circ f_t + iv_t \circ f_t) + o(\varepsilon)$.

To uncouple equations (19) and (20), represent $(u_t, v_t)$ by two functions $\Phi_t$ and $\Psi_t$ using $u_t = \partial_y \Phi_t + \partial_x \Psi_t$ and $v_t = \partial_x \Phi_t - \partial_y \Psi_t$. Notice that $\Phi_t$ is the potential function and $\Psi_t$ is the stream function for the vector field $(v_t, u_t)$. Equations (19) and (20) then become the Poisson equations

$$(21) \qquad \qquad \Delta \Psi_t = 2 \operatorname{Re}(\sigma_t),$$

$$(22) \qquad \qquad \Delta \Phi_t = 2 \operatorname{Im}(\sigma_t).$$

To get a unique solution for these equations, one must specify boundary conditions on $\Phi_t$ and $\Psi_t$. These boundary conditions correspond to different vector field embeddings under the dilatation constraint $\mu_t$.

The examples in the previous sections set $\mu_t := t\mu$ for $t \in [0,1]$ and the boundary conditions $\Phi_t = \Psi_t = 0$ on $\partial \Omega_t$, where $\Omega_t = f_t(\Omega)$. In our implementation of this algorithm, at each time step $t$, we interpolated $\sigma_t$ to a grid enclosing $\Omega_t$, then used a fast Poisson equation solver for grid data with MATLAB.

**6. Discussion.** Our main motivation for developing this methodology is to demonstrate that estimating deformations of arbitrarily smooth isotropic random fields is not only theoretically possible, but can be accomplished in practice. We make no claim to optimality of our methods, but we do hope that this will serve as a test case in the pursuit of efficient and robust methodology for these models in spatial statistics. The use of local likelihood techniques gives quick and easy estimates of $(\mu, \phi)$ that seem to work quite well in simulation. Unfortunately they are not very amenable to theoretical study. We do think, however, that quasiconformal theory presents a promising avenue for the theoretical study and quantification of variability for estimates of deformations. In particular, we believe that quasiconformal theory will be useful for constructing consistent estimates of $C^1$ diffeomorphisms and can also be used for developing flexible parametric models of $C^1$ diffeomorphisms.

We conclude this section with some simulations. First, we demonstrate our methods on differentiable processes. Figure 7 shows the results of applying our methodology to a simulation of a deformed differentiable isotropic process. The true deformation is generated by a vector field flow and the original isotropic process has autocovariance function $K$ of the form $K(t) = 0.0231 - (0.4034)t^2 + |t|^3 + o(|t|^3)$ as $|t| \to 0$ and is simulated on a $400 \times 400$ square grid on $[0,1]^2$; see [2] for details. Notice that even though the process is smooth and, thus, it is more difficult to see the deformation, there is sufficient information embedded in higher order increments for very accurate estimation of the deformation.

Another aspect of these models that warrants further study is the construction of distances on deformations that would allow comparison of methodologies and a study of the variability of the deformation estimates. Since the deformations are only estimable up to rigid motions of the plane, we want these measures to be invariant under such postcomposition. Letting $\Omega$ denote the observation region in $\mathbb{C}$, the following are two candidates that have this invariance:

$$d_1^2(\hat{f}, f) = \int\int_{\Omega^2} (|\hat{f}(z) - \hat{f}(w)| - |f(z) - f(w)|)^2 \, dz \, dw,$$

$$d_2^2(\hat{f}, f) = \int_\Omega |\mu_{\hat{f}}(z) - \mu_f(z)|^2 \, dz.$$

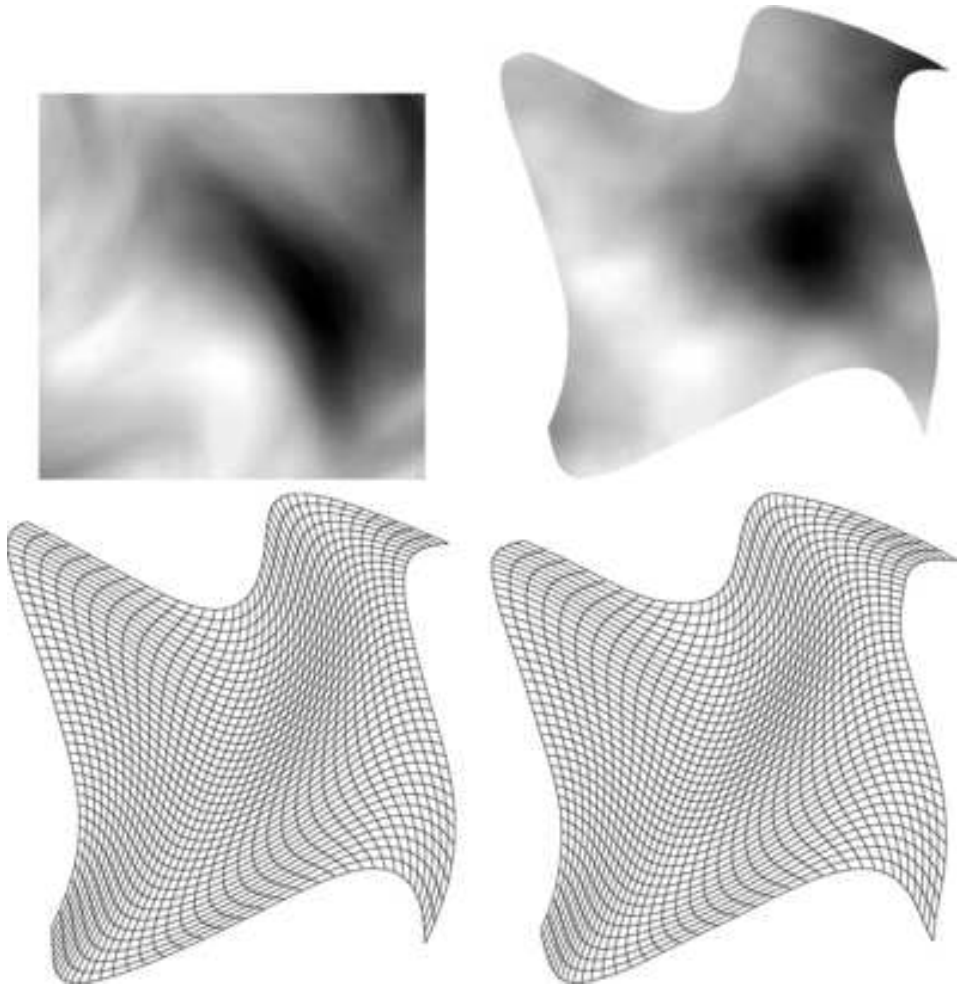

FIG. 7. *Deformation estimation for differentiable processes.* Upper left: *The deformed process $Z \circ f^{-1}$.* Upper right: *The estimate of the original isotropic process $Z$ using the estimated deformation $\widehat{f^{-1}}$.* Bottom left: *The true deformation $f^{-1}$.* Bottom right: *The estimated deformation $\widehat{f^{-1}}$.*

The first distance is based on how well the interpoint distances of $\hat{f}$ match $f$. The second distance is a function of the dilatations of $\hat{f}$ and $f$, instead of the pointwise values.

As a first step toward quantifying the variability of the deformation estimates, we independently simulated the same deformed process as in Section 4 three times. We then added white noise to the last two simulations with a standard deviation of 10% and 25% of the process standard deviation, respectively. The simulations and the corresponding deformation estimates are shown in Figure 8. First notice that the estimated deformation in Figure 5 and the estimated deformation shown top right in Figure 8 are from the same estimation procedure applied to two independent samples. This gives a sense for the variability of the methodology applied to this deformed process. Notice the large scale structure of the deformation looks to be qualitatively stable. Also, notice that one could improve the deformation estimates from the last two samples by modeling the white noise in the likelihood estimation steps. We find it interesting that even without modeling the white noise the estimates are reasonable.

TABLE 1
*Two different distance measures on deformations (columns). The rows correspond to the simulations and estimates shown in Figure 8*

|             | $d_1(\widehat{f^{-1}}, f^{-1})$ | $d_2(\widehat{f^{-1}}, f^{-1})$ |
|-------------|-------------------------------|-------------------------------|
| No noise    | 0.0563                        | 0.0675                        |
| 10% noise   | 0.0655                        | 0.0842                        |
| 25% noise   | 0.1257                        | 0.1354                        |

For each estimate, we also computed its distance from the truth for each of our proposed distance measures (the integrals being approximated by Riemann sums). The results are displayed in Table 1. Notice that adding white noise degrades the deformation estimates. This can be seen both qualitatively and from the distance measurements. We think, however, that the noisy estimates are reasonable and show a certain robustness and stability.

## APPENDIX:
## QUASICONFORMAL MAPS AND COMPLEX DILATATIONS

The goal of this section is to give a short introduction to the theory of quasiconformal maps and to state, without proof, the Mapping theorem (Ahlfors [1], Chapter 5) concerning the uniqueness of a map with a given dilatation. A complete study of quasiconformal maps can be found in Ahlfors [1], Lehto and Virtanen [15], Krushkal' [13] or Ławrynowicz [14]. The Mapping theorem will show that the complex dilatation, $\mu$ of an orientation preserving $C^1$ diffeomorphism, $f$, with bounded distortion uniquely determines $f$ up to postcomposition with a conformal map. The Mapping theorem is more general, however, and will not, for example, require the maps to be differentiable everywhere as in the case of $C^1$ diffeomorphisms. We first define $C^1$ diffeomorphisms and take advantage of the extra smoothness assumptions to redefine $\mu$. We then derive some properties of the complex dilatation and show how $\mu$ characterizes the infinitesimal ellipse that gets mapped to an infinitesimal circle under $f$. Finally, we will relax the $C^1$ diffeomorphic assumptions to state the Mapping theorem.

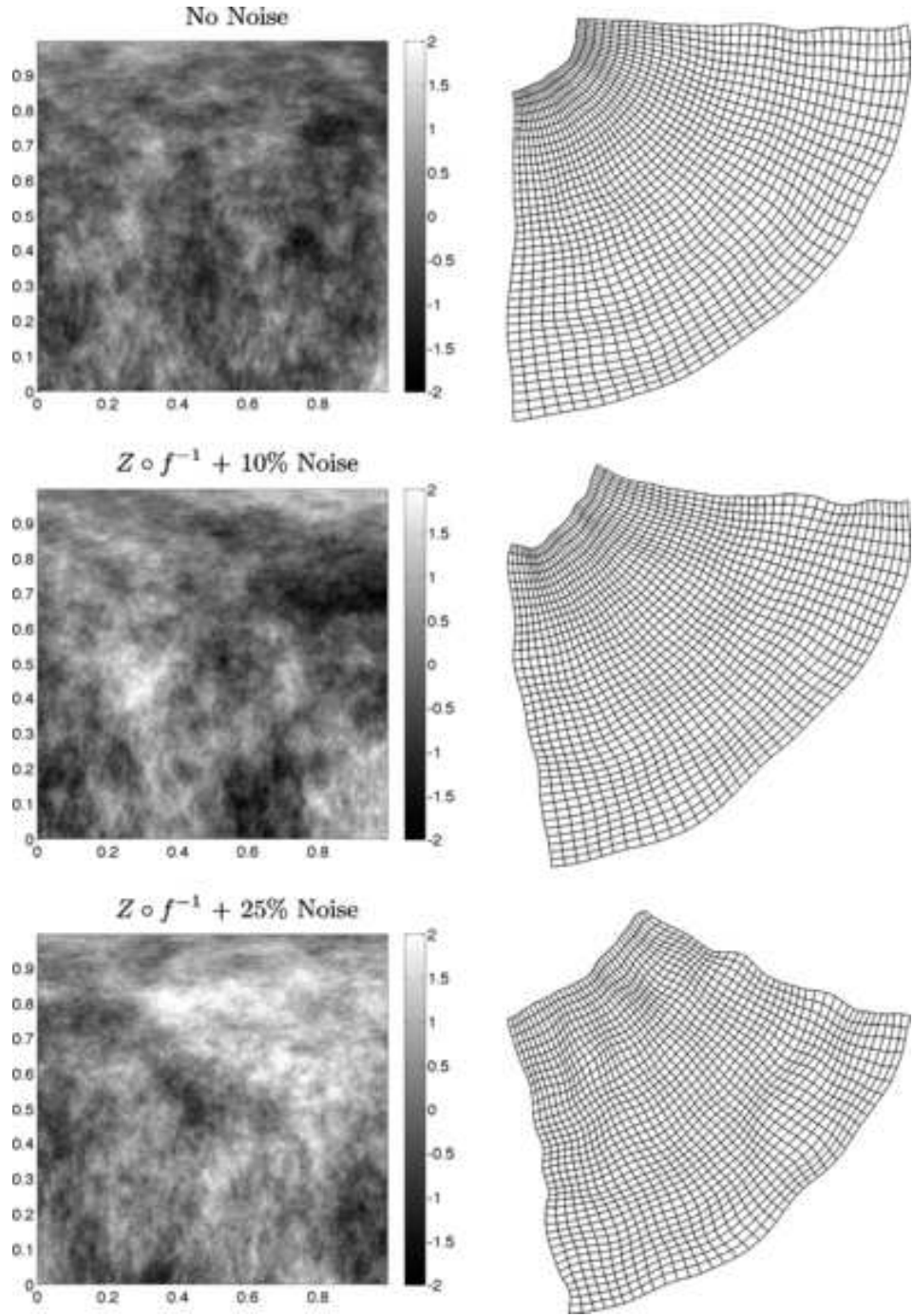

FIG. 8. *Independent simulations of the deformed process from Section 4 with white noise (left column) and the corresponding estimated deformation (right column).*

We will be considering maps $f$ that are defined, not only on $\mathbb{C}$, but also on domains $\Omega \subset \mathbb{C}$. We choose to consider only well-behaved domains, which will allow succinct exposition. In particular, by a *region*, we mean a subset $\Omega \subset \mathbb{C}$ that is either the whole plane $\mathbb{C}$ or is a simply connected open subset of $\mathbb{C}$ whose boundary, $\partial \Omega$, is a Jordan closed curve (also known as a Jordon region, see [20], e.g.). Define a homeomorphism $f$ from a region $\Omega$ to $\Omega'$ to be a continuous, one-to-one, onto map from $\Omega$ to $\Omega'$ such that both $f$ and $f^{-1}$ are continuous. All the maps we wish to consider are necessarily homeomorphic; however, homeomorphisms are not smooth enough to allow us to talk about directional derivatives and Jacobians, which provide the easiest way to understand a complex dilatation. For this reason, we first restrict our homeomorphisms to be sufficiently smooth so that derivatives are defined everywhere and are continuous, that is, $C^1$ diffeomorphic. A homeomorphism $f(\mathbf{x}) = (u(\mathbf{x}), v(\mathbf{x}))$ is an orientation preserving $C^1$ diffeomorphism if $u$ and $v$ are continuously differentiable and for any $\mathbf{x} \in \Omega$,

$$(23) \qquad f(\mathbf{x} + \mathbf{h}) = f(\mathbf{x}) + J_f \mathbf{h} + o(|\mathbf{h}|) \qquad \text{as } |\mathbf{h}| \to 0,$$

where $J_f$ is the Jacobian of the map $f$ and $\det J_f > 0$. Now since all the derivatives exist, we can define $\partial_z f$ and $\partial_{\bar{z}} f$ as in Section 5. Notice that $\partial_z f \neq 0$ since $\det J_f > 0$ and $\det J_f = u_x v_y - u_y v_x = |\partial_z f|^2 - |\partial_{\bar{z}} f|^2$. This allows us to define the complex dilatation $\mu = \mu_f := \partial_{\bar{z}} f / \partial_z f$. We will see that $\mu(z)$ agrees with the complex dilatation of the map defined in Section 3 as parameterizing the eccentricity and inclination of the infinitesimal ellipse centered at $z$ that gets mapped to an infinitesimal circle under $f$.

To see how $\mu$ characterizes the local behavior of the diffeomorphism $f$, we switch to complex notation, $z = x + iy$, so we can write the behavior of $f$ as

$$f(z) = f(z_0) + f_x(z_0)(x - x_0) + f_y(z_0)(y - y_0) + o(|z - z_0|).$$

Rearranging terms, we get $f(z) = f(z_0) + \partial_z f(z_0)(z - z_0) + \partial_{\bar{z}} f(z_0)(\overline{z - z_0}) + o(|z - z_0|)$ and then factoring the nonzero term $\partial_z f$, we can decompose the local behavior of $f$ into an initial stretch map, a final rotation and a uniform stretching

$$(24) \qquad f(z_0 + h) = f(z_0) + \partial_z f(z_0)(h + \mu \bar{h}) + o(|h|),$$

where $\mu = \partial_{\bar{z}} f / \partial_z f$ is the complex dilatation. Now define the directional derivative of $f$ at an angle $\beta \in [0, 2\pi)$ by

$$\partial_\beta f(z) = \lim_{\epsilon \to 0} \frac{f(z + \epsilon e^{i\beta}) - f(z)}{\epsilon e^{i\beta}},$$

the limit existing by the differentiability of $f$. Notice that when $\mu$ is zero at a point $z_0$, (24) tells us that the directional derivative of $f$ at $z_0$ does not depend on direction, that is, that $f$ is conformal at $z_0$. In particular, $\partial_{\bar{z}} f = 0$ reduces to the Cauchy–Riemann equations. Now since $\partial_\beta f(z) = e^{-i\beta} \frac{\partial}{\partial \epsilon} f(z + \epsilon e^{i\beta})|_{\epsilon=0}$, we get $\partial_\beta f = e^{-i\beta}(f_x \cos\beta + f_y \sin\beta) = \partial_z f + e^{-2i\beta} \partial_{\bar{z}} f$. Therefore, the magnitudes of maximal and minimal stretching are attained at $\max_\beta |\partial_\beta f| = |\partial_z f| + |\partial_{\bar{z}} f|$ and $\min_\beta |\partial_\beta f| = |\partial_z f| - |\partial_{\bar{z}} f|$, where the last equality holds since $\det J_f = |\partial_z f|^2 - |\partial_{\bar{z}} f|^2 > 0$ implies that $|\partial_z f| > |\partial_{\bar{z}} f|$. The singular value decomposition of $J_f$ yields the representation (8) with $\lambda_1 \geq \lambda_2 > 0$, which equal the major and minor magnitudes of local stretching. Therefore,

$$\frac{\lambda_1}{\lambda_2} = \frac{\max_\beta |\partial_\beta f|}{\min_\beta |\partial_\beta f|} = \frac{1 + |\mu|}{1 - |\mu|}.$$

To see how $\mu$ is related to the inclination of the ellipse, notice that $\partial_\beta f = \partial_z f(1 + \mu e^{-2i\beta})$ so that $\max_\beta |\partial_\beta f|$ is attained when $2\beta = \arg(\mu)$. Therefore, the inclination of the ellipse that gets mapped to the circle is $\arg(-\mu)/2$, which agrees with (9).

In Section 3 we defined a measure of distortion at a point $z_0$, induced by $f$, to be the ratio of $\limsup_{z \to z_0} |f(z) - f(z_0)|/|z - z_0|$ to $\liminf_{z \to z_0} |f(z) - f(z_0)|/|z - z_0|$. Since this ratio is equal to $\max_\beta |\partial_\beta f|/\min_\beta |\partial_\beta f|$, our original notion of bounded distortion at a point $z_0$ is equivalent to $|\mu(z_0)| < 1$. Therefore, the assumption that $f$ has uniformly bounded distortion amounts to the supposition that $\|\mu_f\|_\infty < 1$.

It turns out that our smoothness assumptions on $f$ are unnatural and the presentation of the theory is best done in full generality. Instead of forcing the derivatives $u_x$, $u_y$, $v_x$ and $v_y$ to exist everywhere, we only suppose existence in the distributional sense and that they be locally integrable (see [8]). Here is the full definition of a quasiconformal map taken from [8]:

DEFINITION A.1. *An orientation preserving homeomorphism $f$ from a region $\Omega$ to a region $f(\Omega)$ is quasiconformal if there exists $k < 1$ such that $f$ has locally integrable, distributional derivatives $\partial_z f$ and $\partial_{\bar{z}} f$ on $\Omega$, and $|\partial_{\bar{z}} f| \leq k |\partial_z f|$ almost everywhere.*

In particular, any orientation preserving $C^1$ diffeomorphism with uniformly bounded distortion is by definition considered a quasiconformal map. Now we have the following existence and uniqueness theorem, taken from [15], page 194.

THEOREM A.1. *Let $\Omega$ and $\Omega'$ be conformally equivalent regions and $\mu$ a measurable function in $\Omega$ with $\|\mu(z)\|_\infty < 1$. Then there exists a quasiconformal mapping $f : \Omega \to \Omega'$ whose complex dilatation coincides with $\mu$ almost everywhere. This mapping is uniquely determined up to a conformal mapping of $\Omega'$ onto itself.*

DEFORMATIONS AND RANDOM FIELDS 23

A useful restatement of the above theorem is that all quasiconformal maps on $\Omega$ can be represented by finding a map $f$, with dilatation $\mu$, then postcomposing it with conformal maps.

We conclude by mentioning that for a general dilatation $\mu$ the associated quasiconformal maps will always be a homeomorphism, but not generally a $C^1$ diffeomorphism. There are sufficient conditions on the complex dilatation $\mu$ to guarantee $C^1$ diffeomorphic solutions; see [15], Theorem 7.2, for example.

DEPARTMENT OF STATISTICS
UNIVERSITY OF CALIFORNIA, BERKELEY
367 EVANS HALL
BERKELEY, CALIFORNIA 94720-3860
USA
E-MAIL: anderes@stat.berkeley.edu

DEPARTMENT OF STATISTICS
UNIVERSITY OF CHICAGO
5734 UNIVERSITY AVE.
CHICAGO, ILLINOIS 60637
USA
E-MAIL: stein@galton.uchicago.edu